%% file: elem4.tex
\documentclass[a4paper,11pt]{amsart}
\usepackage{amsmath,amsthm,amssymb}
\usepackage{color}
\usepackage{graphicx}
\DeclareMathOperator{\Rea}{Re}
\newcommand{\anb}{\bo\partial_{\mathrm{an}}}
\DeclareMathOperator{\sgn}{sgn}

\DeclareMathOperator{\Arg}{Arg}
\newcommand{\bo}{\boldsymbol}
\newcommand{\qand}{\qquad\text{and}\qquad}
\theoremstyle{definition}
\newtheorem{point}{}[section]

\theoremstyle{plain}

\newtheorem{lemma}[point]{Lemma}
\newtheorem{theorem}[point]{Theorem}
\newcommand{\marginextend}[1]{ \addtolength{\oddsidemargin}{-#1}  \addtolength{\evensidemargin}{-#1}
  \addtolength{\textwidth}{#1}\addtolength{\textwidth}{#1}}
\newcommand{\updownextend}[1]{ \addtolength{\topmargin}{-#1}  \addtolength{\textheight}{#1}
\addtolength{\textheight}{#1}}
\marginextend{1cm}
\updownextend{0cm}
\title{Spectral calculations on locally convex vector spaces I}
\author{Gyula Lakos}
\address{Department of Geometry, E\"otv\"os University, P\'azm\'any P\'eter s.~1/C,  Budapest, H--1117, Hungary}
\email{lakos@cs.elte.hu}
\keywords{Functional calculus, locally convex vector spaces, Lipschitz  curves.}
\subjclass[2000]{Primary: 47A60, Secondary: 47A06.}
\begin{document}
\begin{abstract}
We develop a holomorphic functional calculus for multivalued linear operators on locally convex vector spaces,
based on the resolvent identity.
This includes the case of fractional powers along Lipschitz   curves.
\end{abstract}
\maketitle

\section*{Introduction}
The objective of  this paper is to  elaborate the following construction of a functional calculus:
If $\Gamma$ is a Lipschitz  curve system on the Riemann sphere $\overline{\mathbb C}$, then any
slowly growing holomorphic function on $\overline{\mathbb C}\setminus \Gamma$ can be represented
as the integral of a ``resolvent measure'' on $\Gamma$. In case of $\overline\Gamma\subset\mathbb C$, it
looks like
\[f(z)=\int_{w\in\Gamma}\sum_{i=0}^n\frac{\mu_i(w)}{(w-z)^i},\]
where $n\in\mathbb N$, and $\mu_i$ are appropriate complex Borel measures on $\Gamma$ of finite variation.
Now, if  $A$ is a not necessarily everywhere defined operator on a locally convex vector space $\mathfrak V$,
and $u\in \mathfrak V$, then it is reasonable to define $f(A)u$ as the action
of the integral of the resolvent measure on $u$, where $A$ is substituted to the place of $z$.
 In case of $\overline\Gamma\subset\mathbb C$, this yields
\[f(A)u=\int_{w\in\Gamma}\sum_{i=0}^n\frac{\mu_i(w)}{(w-A)^i}u.\]
In particular, our functional calculus is based on the resolvent identity and
related complex function theory, in contrast to other kinds of approaches based, for example,
on power series expansions or monotonicity.

In a more detailed way:
The slowly growing holomorphic functions on the complement of the Lipschitz   curve system
$\Gamma\subset\overline{\mathbb C}$
form the algebra $\mathcal H_\Gamma$.
This algebra has a natural filtration: $f\in\mathcal H_\Gamma^n$, if $f$ is an appropriate
superposition of rational functions with poles in $\Gamma$ with multiplicity at most $n$.
If $A$ is a not necessarily everywhere defined operator on a locally convex vector space $\mathfrak V$,
 then $u\in\mathrm D^n_\Gamma(A)$ if $u\in\mathfrak V$ and
the resolvent operators of $A$ with poles from $\Gamma$ applied at most $n$ times to $u$ give
appropriately bounded continuous functions on $\Gamma^n$.
Then, we can show that in case of $f\in\mathcal H_\Gamma^n$, $u\in\mathrm D^n_\Gamma(A)$, the expression
$f(A)u$ is well-defined. Furthermore, we demonstrate the
viability of this functional calculus in arithmetical sense, which
includes the statements of  global conformal invariance, linearity properties, and multiplicativity.
Locality, which is essentially the statement that if $f$ vanishes on the spectrum of $u$ with respect to $A$,
then $f(A)u=0$, is relatively technical, hence it will be dealt in  a subsequent paper.

Inspiration for this work came from three sources:
one is the author's interest in taking square roots along special curves in locally convex algebras,
cf.~\cite{L1};  another one is the classical theory of fractional powers of
operators with radial spectral cuts, see  Mart\'\i nez Carracedo and Sanz Alix \cite{MS};
and the third one is the classical Riesz functional calculus, cf.~Dunford and Schwartz \cite{DS}.
For a comprehensive review of related matters, and much of the history of functional calculus,
see the book of Haase \cite{H2},
even if that one concentrates on the sectorial case.

The author would like to thank Markus Haase for some useful comments.

\section{Lipschitz   curves and curve systems}
\begin{point}[Definitions and basics]
A 1-dimensional compact simplicial complex $\mathcal C$ is just finitely many
closed intervals possibly glued together at endpoints, plus finitely many further points.
One can consider a metric space structure on $\mathcal C$ induced from the
natural metric of the intervals.
This metric depends on choices but it is unique up to bi-Lipschitz  equivalence.
We can just take a choice $d_{\mathcal C}$.
We call this $(\mathcal C, d_{\mathcal C})$  a compact 1-dimensional metric simplicial complex.
We can also consider the  Riemann sphere $\overline{\mathbb C}$.
This space can be endowed by a spherical metric $d_{\overline{\mathbb C}}$
in a natural way. This metric is unique only up to conformal equivalence;
nevertheless, in particular, it is unique up to bi-Lipschitz  equivalence.
A compact Lipschitz  curve system $\gamma$
is  a bi-Lipschitz  embedding
$\gamma:(\mathcal C, d_{\mathcal C})\rightarrow(\overline{\mathbb C},d_{\overline{\mathbb C}})$
of a compact 1-dimensional metric simplicial complex into the Riemann
sphere. So, there exist constants $C_1,C_2>0$ such that
$C_1 d_{\mathcal C}(t_1,t_2)\leq  d_{\overline{\mathbb C}}(\gamma(t_2),\gamma(t_1))\leq C_2d_{\mathcal C}(t_1,t_2)$ holds for all $t_1,t_2\in\mathcal C$.
If $\Gamma=\gamma(\mathcal C)$, then it it easy to see that
the metric $d_{\overline{\mathbb C}}$ restricted to $\Gamma$,
the inner (arc length) metric  of $\Gamma$ induced from $d_{\overline{\mathbb C}}$
(extended appropriately from the connected components),
and the  metric $\gamma_*d_{\mathcal C} $ are all  bi-Lipschitz  equivalent to each other.
This means that we can obtain $\gamma$ from $\Gamma$ up to bi-Lipschitz  reparametrization.
So, it is just a matter of convenience to use $\gamma$ or $\Gamma$.

A slight generalization is as follows:  A finite 1-dimensional metric simplicial complex
$(\mathcal F, d_{\mathcal F})$ is just a compact 1-dimensional metric simplicial complex
but with finitely many points omitted. One can see that any such $(\mathcal F, d_{\mathcal F})$
can be completed to a compact 1-dimensional metric simplicial complex
$(\mathcal C, d_{\mathcal C})\equiv (\overline{\mathcal F}, d_{\overline{\mathcal F}})$
uniquely, in which case the omitted points  from $(\mathcal C, d_{\mathcal C})$ are non-isolated.
A finite Lipschitz  curve system $\gamma$
is  a bi-Lipschitz  embedding $\gamma:(\mathcal F, d_{\mathcal F})
\rightarrow(\overline{\mathbb C},d_{\overline{\mathbb C}})$ which, however, extends to
a  bi-Lipschitz  embedding $\overline\gamma:(\overline{\mathcal F},
d_{\overline{\mathcal F}})\rightarrow(\overline{\mathbb C},d_{\overline{\mathbb C}})$.
Of course, the principal examples
for Lipschitz  curve systems are  the (bi-) Lipschitz  curves
$\gamma:[a,b]\rightarrow \overline{\mathbb C}$,  $\gamma:(a,b]\rightarrow \overline{\mathbb C}$,
 $\gamma:\mathbb S^1\rightarrow \overline{\mathbb C}$.
In fact, Lipschitz  curve systems can be thought as collections of Lipschitz  curves  which
join appropriately at certain nodes.

If $\overline{\Gamma}\subset \mathbb C$, then the metric
$d_{\overline{\mathbb C}}$ can be replaced by the standard metric $d_{\mathbb C}$.
If  $\overline\Gamma\subset\mathbb C$ is connected, $0<C$, then we say that
$\Gamma$ has wiggling constant $C$ if for any $s_1,s_2\in\Gamma$ the equality
$ Cd_{\mathbb C,\overline\Gamma} (s_1,s_2)\leq d_{\mathbb C} (s_1,s_2)$ holds, where
$d_{\mathbb C,\overline\Gamma}$ denotes the inner arc-length distance in $\overline\Gamma$.
Lipschitz curves are almost everywhere differentiable. More precisely, if $\gamma:[a,b]\rightarrow\mathbb C$
is a Lipschitz curve such that $C_1|t_1-t_2|\leq  d_{\mathbb C}(\gamma(t_2),\gamma(t_1))\leq
C_2|t_1-t_2|$ holds, then $\gamma'$ exists almost everywhere, and  $C_1\leq |\gamma'|\leq C_2$.
In particular, the differentiable points are dense in $\Gamma=\gamma([a,b])$.

A point of $\Gamma$ is called a (real) node if the number of legs joining to it is not $2$.
(In particular, isolated points of $\Gamma$ are nodes.)
Let $N(\Gamma)$ denote the set of nodes of $\Gamma$; this is a finite set.
We define the set of empty nodes of $\Gamma$ as $\overline{\Gamma}\setminus\Gamma$;
this is also a finite set.  Points of $\Gamma\setminus N(\Gamma)$ are the regular
points of $\Gamma$, they lie on open Lipschitz curves joining two real or empty nodes, or on closed Lipschitz curves.
In what follows, $P_2(\Gamma)$ denotes the subsets of cardinality $2$ of $\Gamma$,
and $P_{\mathrm L}(\Gamma)$ denotes the  set of those subsets  of $\Gamma$
which are Lipschitz curves themselves. Both sets are easy to topologize; in fact, they have nice cellular structures.
\end{point}
\begin{point}[The analytic boundary]
Suppose that $\Gamma$ is a compact Lipschitz curve system, and $\mathcal U$ is the union of some
connected components of $\overline{\mathbb C}\setminus\Gamma$. Then $\mathcal U$ is an open set.
If $\mathcal U=\overline{\mathbb C}\setminus\Gamma$, then its
boundary is $\bo\partial \mathcal U=\Gamma$, and in general,  $\bo\partial \mathcal U\subset\Gamma$.
However, there is a richer boundary, the analytic boundary $\anb \mathcal U$.
A sequence $x_n$ in $\mathcal U$
converges to a point of the analytic boundary, if
it converges to a point $x$ of the boundary,
and it eventually stays in a single local connected component of $\mathcal U$ around $x$.
We can topologize  $\anb \mathcal U$,  so  $\mathcal U\cup\anb \mathcal U$
is  a topological manifold with boundary, and $\anb \mathcal U$
is  the topological boundary plus finitely many points.
There is a surjective continuous forgetful map
$\mathsf f:\anb \mathcal U\rightarrow \bo\partial\mathcal U$,
which is typically not injective.
If $\mathcal U=\overline{\mathbb C}\setminus\Gamma$ and $s\in\Gamma$ has $n\geq 1$ legs joining to $s$, then
$\mathsf f^{-1}(s)$ has exactly $n$ points;  if $s\in\Gamma$ has $n=0$ legs joining to $s$,
i.~e.~it is an isolated point, then $\mathsf f^{-1}(s)$ still has  $1$ point.
If $\mathcal U$ is a connected component of $\overline{\mathbb C}\setminus\Gamma$,
then we say that it is a Lipschitz domain.
In general, the analytic boundary of $\mathcal U$ is  the disjoint union of the
analytic boundary of its connected components.
Topologically, each component of $\anb \mathcal U$ is either a point or a circle.
In fact, for each connected component $K_i$ of $\anb \mathcal U$
we can find an closed, oriented boundary curve $\hat\gamma_i:\mathbb S^1\rightarrow K_i$, such that
$\hat\gamma_i$ degenerates to a point, or the underlying curve
$\gamma_i=\mathsf f\circ\hat\gamma_i$ is a piecewise Lipschitz   curve, although not necessarily simple,
and $\mathcal U$ ``lies  on the left'' of $\hat\gamma_i$; in particular,
in case of a bounded, simply connected $\mathcal U$, it goes counterclockwise.
These oriented boundary curves are unique up to bi-Lipschitz  reparametrizations.
If $\mathcal U$ is an $n$-ply connected Lipschitz domain, then  $\anb \mathcal U$
has exactly $n$ connected components.

 The simplest example is when $\gamma:[a,b]\rightarrow\mathbb C$ is a Lipschitz curve,
$\Gamma=\gamma([a,b])$.
Going counterclockwise, the analytic boundary  of $\mathbb C\setminus\Gamma$ in $\mathbb C$ contains
the endpoint $\gamma(b)$, the left boundary $\Gamma^+=\{\gamma(t)^+\,:\,t\in(a,b)\}$,
the endpoint $\gamma(a)$, and the right boundary $\Gamma^-=\{\gamma(t)^-\,:\,t\in(a,b)\}$.
The same notation applies when $s$ is a regular point of a possibly even non-compact $\Gamma$, and $\Gamma$
is oriented in a neighborhood $\mathcal U$ of $s$.
Then  $s^+$ denotes the point of the analytic boundary of $\overline{\mathbb C}\setminus\Gamma$
which lies over $s$, and it is on left side of the curve. Similarly, $s^-$ lies on the right side.
(In general, for non-compact $\Gamma$, the status of the empty nodes is
ambiguous with respect to the analytic boundary.)

If $f$ is a continuous function on $\mathcal U$, then we might want to extend it to the analytic boundary
$\anb(\mathcal U\setminus\Gamma)|_{\mathcal U}\equiv \anb(\mathcal U\setminus\Gamma)
\cap\mathsf f^{-1}(\mathcal U)$.
We can do this according to the notion of convergence we discussed.
The situation is the nicest when $f$ extends continuously to open subsets of the analytic boundary.
Then, for the extended values, the notation  $f(z^+)$, and $f(z^-)$ is reasonable, whenever applicable.
\end{point}
\begin{point}[Loop structure] Let $b_0(\Gamma)$ denote the number of
connected components of $\Gamma$, and let $b_1(\Gamma)$ denote the number of
(homologically) independent loops in $\Gamma$.
Then $b_1(\Gamma)+1$ is the number of connected components of $\overline{\mathbb C}\setminus\Gamma$.
Suppose that $\Gamma\subset\mathbb C$ is compact, and
$\mathcal U_k$, $k=0,\ldots,b_1(\Gamma)$ are the connected components of
 $\overline{\mathbb C}\setminus\Gamma$. Then we can take $\anb \mathcal U_k$
with its canonical orientation (on the non-degenerate pieces),
and we can take their arc-variation measure $\mathrm d \anb \mathcal U_k$
(where degenerate boundary components have no contribution). Their sum
$\sum_{k=0}^{b_1(\Gamma)}\mathrm d \anb \mathcal U_k$ is equal to $0$ 
(because every regular point $\Gamma$ is sided by
two curves of opposite orientation), but any $b_1(\Gamma)$ of those measures represent
independent loops in linear sense. We can assume that $\infty\in\mathcal U_0$,
and then the boundary arc-variation measures $\mathrm d \anb \mathcal U_k$
($k=1,\ldots,b_1(\Gamma)$)
of the bounded connected components of $\overline{\mathbb C}\setminus\Gamma$ represent the
independent loops of $\Gamma$.
The situation is similar in the case when $\Gamma$ is non-compact, except then
for a connected component $\mathcal V$ of $\overline{\mathbb C}\setminus\Gamma$,
the set $\mathcal V\setminus (\overline\Gamma\setminus\Gamma)$ is a disjoint union
$\mathcal U_{\mathcal V,1}\cup\ldots\cup \mathcal U_{\mathcal V,r}$ of some connected components of
$\overline{\mathbb C}\setminus\overline\Gamma$; and we set $ \mathrm d \anb \mathcal U=
\mathrm d \anb\mathcal U_{\mathcal V,1}+\ldots+ \mathrm d \anb\mathcal U_{\mathcal V,r}$.
These variation measures always vanish  around the points of $\overline\Gamma\setminus\Gamma$.
\end{point}
\begin{point}[Measures on Lipschitz curve systems]
Suppose that $\Gamma\subset\overline{\mathbb C}$ is a Lipschitz curve system.
Let $\mathrm M_{\Gamma}^n$ denote the space of such measures $\mu$ on $\Gamma$ such that
(i) $\mu $ complex Borel measure of finite variation on $\Gamma$; and
(ii) for all $w\in\overline{\Gamma}\setminus\Gamma$, the measure
$d_{\overline{\mathbb C}}^n(w,s)^{-1}\mu(s)$
is still a complex Borel measure of finite variation.
This is the same thing as to say that
$\left(\prod_{w\in\overline{\Gamma}\setminus\Gamma}d_{\overline{\mathbb C}}^n(w,s)\right)^{-1}\mu(s)$
is a complex Borel measure of finite variation on $\Gamma$.

If $\overline\Gamma\subset\mathbb C$, then this is equivalent to say that
$Z_\Gamma(s)^{-1}\mu(s)$ is a complex Borel measure of finite variation on $\Gamma$, where
$Z_\Gamma(s)=\prod_{w\in\overline{\Gamma}\setminus\Gamma}(s-w)$.
\end{point}
\begin{point}[Resolvent measures]
If $\overline\Gamma\subset\mathbb C$, then a resolvent measure of order $n$ is a formal sum
\begin{equation}
\bo\mu(s,X)=\mu_0(s)+\frac{\mu_1(s)}{s-X}+\ldots+\frac{\mu_n(s)}{(s-X)^n},\label{eq:resmea}
\end{equation}
where $\mu_k\in\mathrm M_\Gamma^k$, and $X$ is a formal variable.
Resolvent measures can  also be considered in general, but then,  for $s\sim\infty$
they should be of shape
\[ \bo\mu(s,X)=\mu_0(s)+\frac{\mu_1(s)}{s^{-1}-X^{-1}}+\ldots+\frac{\mu_n(s)}{(s^{-1}-X^{-1})^n}.\]

The point is that there is a natural resolvent vector bundle over $\overline{\mathbb C}$, such that the fiber over
$s\in\mathbb C$ is spanned by $1$ and $\frac1{s-X}$, and the fiber over $\infty$
is spanned by   $1$ and $X$. (This resolvent bundle has an additional structure in the sense that
global conformal transformations induce an action on the bundle, too.)
This resolvent bundle yields a bundle over $\overline\Gamma$, and in this
terminology, a resolvent measure is a measure of finite variation, which takes values
in a symmetric power of the resolvent bundle over $\Gamma$, with metrical structure rescaled around
$\overline\Gamma\setminus\Gamma$.

In fact, if $h(s,X)$ is a continuous, nowhere scalar section of the resolvent bundle over $s\in\overline\Gamma$,
then any resolvent measure over $\Gamma$ can (uniquely) be written in form
\[\bo\mu(s,X)=\mu_0(s)+\mu_1(s) h(s,X)+\ldots+\mu_n(s) h(s,X)^n,\]
where $\mu_k\in\mathrm M_\Gamma^k$.
\end{point}

\section{The calculus of multivalued linear operators}
\begin{point}[Elementary properties]
A multivalued linear operator $A$ on a vector space $\mathfrak V$ is simply a linear subspace of
$\mathfrak V\times\mathfrak V$.
If $A$ is a multivalued linear operator, and $a,b,c,d\in\mathbb C$, then
we may form the multivalued linear operator
\[\frac{aA+b}{cA+d}=\{(ay+bx,cy+dx)\,:\,(y,x)\in A\}. \]
In particular, this yields a $\mathrm{PGL}(2,\mathbb C)$-action on multivalued linear operators.
In simple cases, we use the colloquial notation
$aA$, $A+b$, $b+A$, $A^{-1}$, $(cA+d)^{-1}$, etc.
We have to keep in mind that the definition is a matrix action,
so $a+\frac{b}{c-A^{-1}}$ is to be imagined as at least three consecutive matrix actions
($A\mapsto A^{-1}$, $B\mapsto \frac{b}{c-B}$, $C\mapsto a+C$) executed, but which
simplify according to the associativity of the action.
If $h(z)=\frac{az+b}{cz+d}$ is the corresponding conformal map on the Riemann sphere,
then we use the short notation $h_{\mathrm f}(A)=\frac{aA+b}{cA+d}$.
The multivalued linear operator $A$ is a (not necessarily everywhere defined) linear operator
if $(y,x)\in A, x=0$ implies $y=0$, i.~e.~its values are well-defined.
In that case, we define $Ax$ as $y$, where $(y,x)\in A$.
\end{point}
\begin{lemma}\label{lem:quasi1}
(a) If $a,b\in\mathbb C$, $a\neq 0$, then
$Au$ exists if and only if
$(aA+b)u$ exists,  and in that case $(aA+b)u=aAu+bu$.

(b) Suppose that $Au,Bu$ exist, $a,b\in\mathbb C$.
Then  $BAu$ exists if and only if  $(B+b)(A+a)u$ exists,  and
in that case $(B+b)(A+a)u=abu+bAu+aBu+BAu$.

(c) If $Bu$ and $B^{-1}v$ exist, then $B^{-1}Bu$ and $BB^{-1}v$ exist, moreover
$B^{-1}Bu=u$, $BB^{-1}v=v$.

(d) If $(B+b)(B^{-1}+a)u$ exists, $a,b\in\mathbb C$, $a\neq 0$, then $Bu$ exists. \qed
\end{lemma}
\begin{lemma}\label{lem:quasi2}
Suppose that $a_1d_1-b_1c_1\neq0$, $a_2d_2-b_2c_2\neq0$,
$\frac{a_1X+b_1}{c_1X+d_1}\frac{a_2X+b_2}{c_2X+d_2}
=\frac{a_3X+b_3}{c_3X+d_3}$.
(So,  $a_3d_3-b_3c_3\neq0$ or $\frac{a_3X+b_3}{c_3X+d_3}\equiv c\in\mathbb C\setminus\{0\}$.)
 If

(i) $ \frac{a_1A+b_1}{c_1A+d_1}u$ and $ \frac{a_2A+b_2}{c_2A+d_2}u$ exist, or

(ii) $\frac{a_1A+b_1}{c_1A+d_1}\frac{a_2A+b_2}{c_2A+d_2}u$ exists,  $a_3d_3-b_3c_3\neq0$,
~\\then
\[\frac{a_1A+b_1}{c_1A+d_1}\frac{a_2A+b_2}{c_2A+d_2}u=\frac{a_3A+b_3}{c_3A+d_3}u=
\frac{a_2A+b_2}{c_2A+d_2}\frac{a_1A+b_1}{c_1A+d_1}u,\]
including the statement that every expression in the equation exists.
\begin{proof}
Set $B=\frac{a_1A+b_1}{c_1A+d_1}$.
Then $\frac{a_2A+b_2}{c_2A+d_2}=pB^{-1}+q$ with $p,q\in\mathbb C$, $p\neq0$.
What we have to prove is that
\[B(pB^{-1}+q)u=(p+qB)u=(pB^{-1}+q)Bu, \]
including that every term is well-defined.

(i) We know that $Bu$ and $ (pB^{-1}+q)u$ exist. In particular, $B^{-1}u$ also exists.
Using Lemma \ref{lem:quasi1}.c, we see that $BB^{-1}u$ and $B^{-1}Bu$ also exist, and they are equal to $u$.
Then Lemma \ref{lem:quasi1}.b implies that the equality above holds.

(ii) In that case $q\neq 0$.
Then, according to Lemma \ref{lem:quasi1}.d, $Bu$ exists.
By this, we have reduced the problem to case (i).
\end{proof}
\end{lemma}

\begin{lemma}[Resolvent identity and commutativity] \label{lem:quasi3}
Suppose that $w_1,w_2\in\mathbb C$, $w_1\neq w_2$.
Then, the identity
\[\frac1{w_1-A}\frac1{w_2-A}u=\frac1{w_2-w_1}\left(\frac1{w_1-A} u-\frac1{w_2-A}u\right)=
\frac1{w_2-A}\frac1{w_1-A}u,\]
holds in the sense that if any of the three terms is well-defined, then the
other two ones are also well-defined, and equality holds.
In short terms: The resolvent identity and commutativity hold, if either side makes sense.
\begin{proof}
Let $B=\frac{w_1-A}{w_2-A}$, $\xi=\frac1{w_2-w_1}$.
Then
\[(w_1-A)^{-1}=\xi(B^{-1}-1),\qand(w_2-A)^{-1}=\xi(1-B).\]
In particular, $Bu$ and $B^{-1}u$ exist, hence $B^{-1}Bu$ and $BB^{-1}u$ exist, and
they are equal to $u$.
We have to prove that
\[\xi(B^{-1}-1) \xi(1-B) u=\xi(\xi(B^{-1}-1)u-\xi(1-B)u)=\xi(1-B)\xi(B^{-1}-1)u,\]
i.~e.
\[(B^{-1}-1)(1-B)u=Bu-2u+B^{-1} u=(1-B)(B^{-1}-1)u.\]

(i) If $(w_1-A)^{-1}u$ and $(w_2-A)^{-1}u$ exist, then
the equality follows from  Lemma \ref{lem:quasi1}.b.

(ii)  If  $(w_1-A)^{-1}(w_2-A)^{-1}u$ exists, then the statement is reduced to part (i) by Lemma \ref{lem:quasi1}.d.
\end{proof}
\end{lemma}
If $n\geq0$, $c\in\mathbb C$, $A$ is a multivalued linear operator on $\mathfrak V$, $u\in\mathfrak V$,
then we  use the notation
$\frac{c}{A^n}u=c\stackrel{n\text{ times}}{\overbrace{\frac1A\ldots\frac1A}}u,$
i.~e.~$\frac{c}{A^n}u$ is $A^{-1}$ iterated $n$ times on $u$
multiplied by the constant $c$.

\begin{lemma}[Higher resolvent identities]\label{lem:qRID}
Suppose that $w_1,w_2\in\mathbb C$, $w_1\neq w_2$. Then

(a) For $n\in\mathbb N$,
\[\frac{1}{w_1-w_2}\left(\frac1{(w_1-A)^{n}}u-\frac1{(w_2-A)^{n}}u\right)=
-\sum_{j=1}^{n}\frac1{(w_1-A)^{j}}\frac1{(w_2-A)^{n-j+1}}u,\]
if either side makes sense.

(b) For $n_1,n_2\in\mathbb N$,
\begin{multline}
\frac1{(w_1-A)^{n_1+1}}\frac1{(w_2-A)^{n_2+1}}u=\\=
\sum_{k=1}^{n_1+1}
\frac{\begin{pmatrix}n_1+n_2+1-k\\n_2\end{pmatrix}}{(w_2-w_1)^{n_2+n_1+2-k}}\frac1{(w_1-A)^{k}}u
+\sum_{j=1}^{n_2+1}
\frac{\begin{pmatrix}n_1+n_2+1-k\\n_1\end{pmatrix}}{(w_1-w_2)^{n_1+n_2+2-k}}\frac1{(w_2-A)^{k}}u,
\label{muli1}
\end{multline}
if either side makes sense.
\begin{proof}
These identities are arithmetical consequences of the resolvent identity.
\end{proof}
\end{lemma}

\begin{point}
Assume  that $\mathfrak V$ is a locally convex vector space.
Suppose that $\Gamma\subset\overline{\mathbb C}$, and $A$ is a multivalued linear operator.
For $n\in\mathbb N$, we define $\mathrm D^n_{\Gamma}(A)\subset\mathfrak V$ as the set of elements
$u$ such that for all $0\leq m\leq n$ the functions
\[R_\infty^m:\quad  (z_1,\ldots,z_m)\in(\Gamma\setminus\{\infty\})^m\mapsto
\frac{ (z_1-A)^{-1}\ldots(z_m-A)^{-1}u}{\prod_{s\in\overline\Gamma\setminus\Gamma\setminus\{\infty\}}  (z_1-s)^{-1}\ldots(z_m-s)^{-1} } \]
and
\[R_0^m:\quad  (z_1,\ldots,z_m)\in (\Gamma\setminus\{0\})^m
\mapsto \frac{(z_1^{-1}-A^{-1})^{-1}\ldots(z_m^{-1}-A^{-1})^{-1} u}
{\prod_{s\in\overline\Gamma\setminus\Gamma\setminus\{0\}}  (z_1^{-1}-s^{-1})^{-1}\ldots(z_m^{-1}-s^{-1})^{-1} }\]
are well-defined and bounded continuous functions (in every seminorm of $\mathfrak V$).
(If $\Gamma$ is compact, then denominator is just $1$.)
Notice that $\mathrm D^1_{\{\infty\}}(A)$ is simply the domain of $A$, if it is an operator.

One might think that $u\in\mathrm D^n_{\Gamma}(A)$
is characterized by a global property with respect to $\Gamma$.
This is not the case, it is a local property with respect to $\overline{\Gamma}$.
Indeed, if the $z_i$ are well-separated, then we can use the
resolvent identity to write the products into additive form,  so existence and continuity is
granted.

It is easy to see that   if $h(s,X)$ is a continuous, nowhere scalar section of the resolvent bundle over $s\in\overline\Gamma$, then  $u\in\mathrm D^n_{\Gamma}(A)$ means that
for all $0\leq m\leq n$ the functions
\[R_h^m:\quad  (z_1,\ldots,z_m)\in\Gamma^m\mapsto
\frac{ h_{\mathrm f}(z_1,A)\ldots h_{\mathrm f}(z_m,A) u}{\prod_{s\in\overline\Gamma\setminus\Gamma}  d_{\overline{\mathbb C}}(z_1,s)^{-1}\ldots d_{\overline{\mathbb C}}(z_m,s)^{-1}}\]
are well-defined and bounded continuous functions. One can topologize
the spaces $\mathrm D^n_{\Gamma}(A)\subset\mathfrak V$ according to these bounded continuous functions.
\end{point}

\begin{point}\label{po:measint}
We remark that if $S$ is a $\mathtt T_2$, $\mathtt M_2$, locally compact topological space,
$v$ is a continuous, bounded $\mathfrak V$ valued function on $S$,
$\mu$ is complex Borel measure of bounded variation
on $S$, then the convergence of the integral
\[\int_{s\in S}v(s)\mu(s) \]
is a matter of sequential convergence in $\mathfrak V$.
This is because the integral can be computed with respect to
a single sequence of divisions of $S$ independently from the integrand.
In particular, if $\mathfrak V$ is sequentially complete, then the integral exists.
So, in what follows, we assume that $\mathfrak V$ is a sequentially complete
locally convex vector space. This is not a very restrictive condition,
$\mathfrak V$ can always be completed.

If, furthermore, if $A$ is a sequentially closed linear operator on $\mathfrak V$, such that $Av$ is still
bounded and continuous, then
\[A\left(\int_{s\in S} v(s)\,\mu(s)\right)=\int_{s\in S} Av(s)\,\mu(s)\]
holds, including the statement that both sides make sense.

It is a trivial observation that if $a,b,c,d\in\mathbb C$, $ad-bc\neq0$ and $A$
is sequentially closed, then $\frac{aA+b}{cA+d}$ is sequentially closed, too.
\end{point}
\begin{lemma}[Global conformal change]
\label{qlem:confchange}
Let $h:\overline{\mathbb C}\rightarrow\overline{\mathbb C}$ be a conformal  map,
$h(z)=\frac{az+b}{cz+d},$ $ad-bc\neq0.$

(a) Suppose that $w,h(w)\in\mathbb C$.
Then
\[\frac1{w-A}=\frac{c}{cw+d}+\frac{ad-bc}{(cw+d)^2}\frac{1}{h(w)-h_{\mathrm f}(A)},\]
if either side makes sense.

(b) Suppose that $\Gamma\subset\mathbb C$ is a Lipschitz   curve system, $\mu$ is a measure on it.
Assume that  $h(\Gamma)\subset\mathbb C$. Suppose that $u\in\mathrm D^n_\Gamma(A)$.
Then
\[\int_{w\in\Gamma}\frac{\mu(w)}{(w-A)^n}u
=\sum_{j=0}^n \int_{w\in\Gamma}\binom{n}{j}\left(\frac{c}{cw+d}\right)^{n-j}
\left(\frac{ad-cd}{(cw+d)^2}\right)^j\frac{\mu(w)}{(h(w)-h_{\mathrm f}(A))^j}u\]
holds. \qed
\end{lemma}
\begin{lemma}[Integral formulas for the resolvent] \label{lem:qRIF}
Assume that  $\gamma:[a,b]\rightarrow\mathbb C$ is a piecewise Lipschitz   curve, and
 $\Gamma=\gamma([a,b])$ is a Lipschitz curve system.

(a)
If $n\geq 1$, $u\in\mathrm D^{n+1}_\Gamma(A)$, then
\begin{equation}
\frac1{(\gamma(b)-A)^{n}}u-\frac1{(\gamma(a)-A)^{n}}u=\int_{t=a}^{b}
\frac{-n\,\mathrm d\gamma(t)}{(\gamma(t)-A)^{n+1}}u.\label{aba1}
\end{equation}

(b) Under similar assumptions, suppose that
\[\sum_{j=1}^{n} \frac{c_j}{(-w)^j}
\sim\sum_{j=1}^{n}\frac{c_j^{(\gamma(a))}}{(\gamma(a)-w)^j}
\sim\sum_{j=1}^{n}\frac{c_j^{(\gamma(b))}}{(\gamma(b)-w)^j}\]
modulo $O(w^{-n-1})$ as $w\rightarrow\infty$. Then
\[\sum_{j=1}^{n} \frac{c_j^{(\gamma(b))}}{(\gamma(b)-A)^j}u-
\sum_{j=1}^{n} \frac{c_j^{(\gamma(a))}}{(\gamma(a)-A)^j}u
=\int_{t=a}^{b}\biggl(\sum_{j=1}^n\binom{n-1}{j-1}c_j\gamma(t)^{n-j}\biggr)
\frac{-n\,\mathrm d\gamma(t)}{(\gamma(t)-A)^{n+1}}u.\]

(c) In particular, if  $n\geq 1$, $u\in\mathrm D^{n+1}_\Gamma(A)$, $\gamma(a)=\gamma(b)$,
then  for any polynomial $P$ with degree $\leq n-1$,
\[\int_{t=a}^b P(\gamma(t))\frac{\mathrm d\gamma(t)}{(\gamma(t)-A)^{n+1}}u=0.\]

(d)
If $n_1,n_2\in\mathbb N$, $u\in\mathrm D^{n_1+n_2+2}_\Gamma(A)$, $\gamma(a)\neq\gamma(b)$, then
\begin{multline}
\frac1{(\gamma(a)-A)^{n_1+1}}\frac1{(\gamma(b)-A)^{n_2+1}}u=\\=\int_{t=a}^{b}
\frac{(n_1+n_2+1)!}{n_1!n_2!}
\frac{(\gamma(t)-\gamma(a))^{n_2}(\gamma(b)-\gamma(t))^{n_1}}{(\gamma(b)-\gamma(a))^{n_1+n_2+1}}
\frac{\mathrm d\gamma(t)}{(\gamma(t)-A)^{n_1+n_2+2}}u.\label{form1}
\end{multline}
\begin{proof}
Part (a) follows from Lemma \ref{lem:qRID}.a. Parts (b), (c), (d) are arithmetical consequences
of part (a).
\end{proof}
\end{lemma}

\begin{point}[Superposition of measures]
Suppose that $T$ and $S$ are $\mathtt T_2$, $\mathtt M_2$, locally compact topological spaces,
and we assign to each point  $s\in S$ a measure $\nu_s$ on $T$.
We say that this assignment is continuous if

(i) the variation map $s\mapsto \|\nu_s\|$ is bounded and continuous;

(ii) the map  $s\mapsto \nu_s$ is weakly continuous on $S$.
~\\ In this case, if $\mu$ is a measure on $S$, then we may form the parametrized measure
\[\tilde\theta(t,s)=\nu_s(t)\mu(s),\]
on $T\times S$, most readily to be defined as a Radon measure.
Integrating in the variable $s$, we obtain the superposition measure
\[\theta(t)=\int_{s\in S} \nu_s(t) \mu(s).\]
In fact, both the parametrized and the superposition measures
can also be taken if $s\mapsto\nu_s$ is not necessarily continuous
but $S$ is a disjoint union of the Borel sets $S_1,\ldots,S_p$, such that on each of them
$s\mapsto\nu_s$ is continuous, and the sets $S_j$ are $\mathtt T_2$, $\mathtt M_2$, locally compact spaces
themselves.

Now, if $v(t)$ is a bounded, continuous, locally convex vector space valued function on $T$, then
\begin{equation}
\int_{s\in S} \left(\int_{t\in T} v(t)\,\nu_s(t)\right)\mu(s)
=\int_{(t,s)\in T\times S} v(t)\,\tilde\theta(t,s) =\int_{t\in T} v(t)\,\theta(t) .
\label{eq:fubini}
\end{equation}

We will consider superpositions only in very special cases, where these superpositions
can also be carried out by more direct means.
\end{point}

\begin{point}[Elementary measures $\Omega$]
Suppose that $\Gamma$ is a Lipschitz  curve with endpoints $s_1,s_2$. Then we define
the complex measure $\Omega_{s_1,\Gamma}^n$ on $\Gamma$ by setting
\[\Omega_{s_1,\Gamma}^n(s)=-n \,\mathrm d\Gamma(s),\]
where $\Gamma$ is oriented in the direction from $s_1$ to $s_2$.
With this notation, \eqref{aba1} reads as
\begin{equation}
\frac1{(s_2-A)^{n}}u
=\frac1{(s_1-A)^{n}}u+\int_{s\in\Gamma}\frac{\Omega_{s_1,\Gamma}^n(s)}{(s-A)^{n+1}}u.
\label{aba2}
\end{equation}
\end{point}
\begin{point}[Additive choice functions]
Now, suppose that $\Delta$ is a connected Lipschitz curve system, $\overline\Delta\subset\mathbb C$.
We say that $\phi:\Delta\rightarrow P_{\mathrm L}(\Delta)$ is an additive choice function with respect to
$w\in\Delta$ if

(i) for any $x\in\Delta$ the choice value $\phi(x)$ is a Lipschitz curve connecting $w$ and $x$.

(ii) $\phi$ is Borel measurable.\\
Such choice functions are, of course, easy to find. If $\Delta$ is a tree,
then $\phi$ is unique.
\end{point}
\begin{theorem}[Additive reduction] \label{thm:addred}
(a)
Suppose that $\mu\in\mathrm M^{n}_\Delta$,  and $\phi$
is an additive choice function with respect to $w\in\Delta$.
Then the superposition measure
\[\nu(s)=\int_{x\in\Delta}\Omega_{w,\phi(x)}^{n}(s)\mu(x)\]
exists, moreover $\nu\in\mathrm M^{n+1}_\Delta$.

(b) Furthermore,
assume that $u\in\mathrm D^{n+1}_\Delta(A)$.
Then
\[\int_{s\in\Delta}\frac{\mu(s)}{(s-A)^n}u=
\left(\int_{s\in\Delta}\mu(s)\right)\frac{1}{(w-A)^n}u+\int_{s\in\Delta}\frac{\nu(s)}{(s-A)^{n+1}}u.\]
\begin{proof} (a) Consider the measures $\tilde\Omega_{w,\Gamma}^{n}(s)=\dfrac{Z_\Delta(x)^n}{Z_\Delta(s)^{n+1}}\Omega_{w,\Gamma}^{n}(s)$
where $\Gamma\in P_{\mathrm L}(\Delta)$ with endpoints $w,x$. These depend continuously on $\Gamma$,
hence it is sufficient to show that their variations are uniformly bounded.

Consider first the special case
when $\lambda:[0,b]\rightarrow\mathbb C$ is a Lipschitz curve, in arc-length parametrization, with wiggling constant $C$; $\Delta=\lambda((0,b])$, $p=\lambda(0)$, $w=\lambda(b)$. Then $0<a<b$, $\Gamma=\lambda([a,b])$, $x=\lambda(a)$ can be assumed,
and for the variation
\[\|\tilde\Omega_{w,\Gamma}^{n}\|\leq\int_{t=a}^b\left\| \frac{a^n}{(Ct)^{n+1}}n\,\mathrm dt\right\|\leq \frac1{C^{n+1}}.\]

Now, in general, we can take $\Delta_{\mathrm{amp}}$, which is $\Delta$ but
pieces of the legs joining to points of $\overline\Delta\setminus\Delta$ amputated, yet so that  $w\in \Delta_{\mathrm{amp}}$.
Now, any Lipschitz curve $\Gamma\subset \Delta$ with endpoints $w,x$ stays in $\Delta_{\mathrm{amp}}$
except if $x$ enters into an amputated leg piece $\Lambda$, when $\Gamma\subset \Delta_{\mathrm{amp}}\cup\Lambda$.
In this case $\|\tilde\Omega_{w,\Gamma}^{n}|_{\Delta_{\mathrm{amp}}}\|$ can be estimated rather trivially, based on
the arc length of $\Delta$, while the estimate of $\|\tilde\Omega_{w,\Gamma}^{n}|_{\Lambda}\|$
is qualitatively the same as the special case above.

(b) This is the superposition of the equalities \eqref{aba2}.
\end{proof}
\end{theorem}
At first sight, this additive superposition measure
might look formidable, but, essentially, $\nu$ is  just $n$ times a primitive function of $\mu$
with starting point $w$, multiplied by the variation measure of $\Delta$, and with some
piecewise constant ambiguities due to the loops in $\Delta$.

\begin{point}[Elementary measures $\Xi$]
If $s_1,s_2\in\mathbb C$, $s_1\neq s_2$, $1\leq k\leq n_1+n_2+1$, then we define the measure
\[\Xi_{s_1,s_2;k}^{(n_1,n_2)}(s)=
\frac{\begin{pmatrix}n_1+n_2+1-k\\n_2\end{pmatrix}}{(s_2-s_1)^{n_1+n_2+2-k}}\,\bo\delta_{s_1}(s)+
\frac{\begin{pmatrix}n_1+n_2+1-k\\n_1\end{pmatrix}}{(s_1-s_2)^{n_1+n_2+2-k}}\,\bo\delta_{s_2}(s);\]
where $\bo\delta_{t}$ the Dirac measure supported at $t$.
Then \eqref{muli1} reads as
\begin{equation}
\frac1{(s_1-A)^{n_1+1}(s_2-A)^{n_2+1}}u=\sum_{k=1}^{n_1+n_2+1}\int_{s\in\{s_1,s_2\}}
\frac{\Xi_{s_1,s_2;k}^{(n_1,n_2)}(s)}{(s-A)^{k}}u.
\label{muli2}
\end{equation}

If $\Gamma$ is a Lipschitz  curve in $\mathbb C$ with endpoints $s_1$ and $s_2$,
then we define the measure $\Xi_\Gamma^{(n_1,n_2)}$ on $\Gamma$ by
\[\Xi_\Gamma^{(n_1,n_2)}(s)=\frac{(n_1+n_2+1)!}{n_1!n_2!}
\frac{(s-s_1)^{n_2}(s_2-s)^{n_1}}{(s_2-s_1)^{n_1+n_2+1}}
\,\mathrm d\Gamma(s),\]
where $\mathrm d\Gamma(s)$ is understood with $\Gamma$ oriented in the direction from $s_1$ to $s_2$.
This is the push-forward of the appropriate expression from $\eqref{form1}$
via $\gamma$. In fact, it is easy  to see that $\Xi_\Gamma^{(n_1,n_2)}$ does not depend
on the ordering of the endpoints of $\Gamma$ but only on  $\Gamma$.
Then, \eqref{form1} reads as
\begin{equation}
\frac1{(s_1-A)^{n_1+1}(s_2-A)^{n_2+1}}u=\int_{s\in\Gamma}
\frac{\Xi_\Gamma^{(n_1,n_2)}(s) }{(s-A)^{n_1+n_2+2}}u. \label{form2}
\end{equation}
It is also easy to see, but notable, that if $\Gamma$ has wiggling constant $C$, then for the variation of the measure
\[\|\Xi_\Gamma^{(n_1,n_2)}\|\leq \frac1{C^{n_1+n_2+1}}.\]
We extend the definition to the degenerate case $\Gamma=\{s_0\}$
by setting
 $\Xi_\gamma^{(n_1,n_2)}(t)=\bo\delta_{s_0}(t)$.
 This is still consistent to \eqref{form2} with $s_1=s_2=s_0$.

For a set $\{s_1,s_2\}$ of cardinality $2$ we define
$\Xi^{(n_1,n_2)}_{\{s_1,s_2\};k}$ as $\Xi^{(n_1,n_2)}_{s_1,s_2;k}$ if $1\leq k\leq n_1+n_2+1$, and as $0$
if $k=n_1+n_2+2$.
For a possibly degenerate Lipschitz curve $\Gamma$ with endpoints $s_1,s_2$, we define
$\Xi^{(n_1,n_2)}_{\Gamma;k}$ as $0$ if $1\leq k\leq n_1+n_2+1$, and as $\Xi^{(n_1,n_2)}_{\Gamma}$
if $k=n_1+n_2+2$. Then, for $X=\{s_1,s_2\}$ or $\Gamma$ with endpoints $s_1,s_2$, the equality
\begin{equation}
\frac1{(s_1-A)^{n_1+1}(s_2-A)^{n_2+1}}u=\sum_{k=1}^{n_1+n_2+2}\int_{s\in X}
\frac{\Xi_{X;k}^{(n_1,n_2)}(s) }{(s-A)^{k}}u\label{colla}
\end{equation}
holds. This collates \eqref{muli2} and \eqref{form2}.
\end{point}

\begin{point}[Multiplicative choice functions] \label{po:multsys}
Suppose that $\Delta$ is Lipschitz curve system, $\overline\Delta\subset\mathbb C$.
Then we can define a choice function
\[\psi:\Delta\times\Delta\rightarrow P_2(\Delta)\cup P_{\mathrm L}(\Delta)\]
such that it satisfies the following properties:

(a) If $\psi(s_1,s_2)\in P_2(\Delta)$, then $\psi(s_1,s_2)=\{s_1,s_2\}$.

(b) If $\psi(s_1,s_2)\in P_L(\Delta)$, then the set of endpoints of $\psi(s_1,s_2)$ is $\{s_1,s_2\}$.

(c) There is a constant $\tilde C_1>0$ such that  for all $w\in\overline\Delta\setminus\Delta$
\[|s_1-s_2|\leq \tilde C_1\min\{|s_1-w|,|s_2-w|\},\]
holds when $\psi(s_1,s_2)\in P_L(\Delta)$.

Furthermore, there is a constant $D>0$ such that $|s_1-s_2|<D$, $\psi(s_1,s_2)\in P_L(\Delta)$
implies that $\psi(s_1,s_2)$ is the Lipschitz curve of minimal length
connecting $s_1,s_2$. (This provides a uniform wiggling constant for $\psi(s_1,s_2)$.)

(d) There is a constant $\tilde C_2$ such that $\psi(s_1,s_2)\in P_2(\Delta)$ implies
\[|s_1-s_2|\geq \tilde C_2\min\{1,|s_1-p|,|s_2-p|\,:\,p\in\overline\Delta\setminus\Delta\}.\]
\end{point}

Such choice functions can be chosen. Indeed, we can take $\Delta_{\mathrm{amp}}$, which is $\Delta$ but
short open pieces of the legs joining to points of $\overline\Delta\setminus\Delta$ amputated.
If $s_1$, $s_2$ are in two different amputated leg pieces $\Lambda_1$, $\Lambda_2$, then $\psi(s_1,s_2)=\{s_1,s_2\}$
can be taken. If $s_1$ is in an amputated leg piece $\Lambda_1$ but $s_2\in\Lambda_1\cup \Delta_{\mathrm{amp}}$,
then we can take the unique connecting Lipschitz curve of $s_1,s_2$ if they are close to each other,
and the discrete choice otherwise. If $s_1,s_2\in\Delta_{\mathrm{amp}}$, then we can take
an  appropriate Lipschitz curve of the discrete choice for
$\psi(s_1,s_2)$ depending on whether $s_1,s_2$ are connected in $\Delta_{\mathrm{amp}}$ or not.

In particular, if $\Delta$ is a compact tree, then for $\psi(s_1,s_2)$, we can always take
the unique Lipschitz curve connecting $s_1$ and $s_2$ in $\Delta$.

\begin{theorem}[Multiplicative reduction] \label{thm:multred}

(a)
Suppose that $\mu_1\in\mathrm M^{n_1+1}_\Delta$, $\mu_2\in\mathrm M^{n_2+1}_\Delta$, and $\psi$
is a multiplicative choice function for $\Delta$.
Then the superposition measures
\[\theta_k(s)=\int_{(s_1,s_2)\in\Delta\times\Delta}\Xi_{\psi(s_1,s_2);k}^{(n_1,n_2)}(s)\mu_1(s_1)\mu_2(s_2),\]
$k=1,\ldots,n_1+n_2+2$, exist. Moreover, $\theta_k\in\mathrm M^{k}_\Delta$.

(b) Furthermore,
assume that $u\in\mathrm D^{n_1+n_2+2}_\Delta(A)$, and $A$ is a sequentially closed multivalued linear operator.
Then
\[\int_{s_1\in\Delta} \frac{\mu_1(s_1)}{(s_1-A)^{n_1+1}} \int_{s_2\in\Delta} \frac{\mu_2(s_2)}{(s_2-A)^{n_2+1}}u
=\sum_{k=1}^{n_1+n_2+2} \int_{s\in\Delta} \frac{\theta_k(s)} {(s-A)^k}u.\]
\begin{proof} (a) It is sufficient to check that the variations of the measures
\[\frac{Z_\Delta(s_1)^{n_1+1}Z_\Delta(s_2)^{n_2+1}}{Z_\Delta(s)^k}\Xi_{\psi(s_1,s_2);k}^{(n_1,n_2)}(s)\]
satisfy  uniform estimates independently from $s_1,s_2$.

First, we consider the special case when $a<0<b$,  $\lambda:[a,b]\rightarrow\mathbb C$ is
a Lipschitz curve with wiggling constant $C$, $\Delta=\lambda([a,0))\cup\lambda((0,b])$, $p=\lambda(0)$.
Then
\begin{multline}
\frac{Z_\Delta(s_1)^{n_1+1}Z_\Delta(s_2)^{n_2+1}}{Z_\Delta(s)^k}\Xi_{\{s_1,s_2\};k}^{(n_1,n_2)}(s)=\\
\qquad\qquad\qquad\qquad\qquad\quad=(s_1-p)^{n_1+1-k}(s_2-p)^{n_2+1}
\frac{\begin{pmatrix}n_1+n_2+1-k\\n_2\end{pmatrix}}{(s_2-s_1)^{n_1+n_2+2-k}}\,\bo\delta_{s_1}(s)+\\+
(s_1-p)^{n_1+1}(s_2-p)^{n_2+1-k}
\frac{\begin{pmatrix}n_1+n_2+1-k\\n_1\end{pmatrix}}{(s_1-s_2)^{n_1+n_2+2-k}}\,\bo\delta_{s_2}(s).
\notag\end{multline}
Now, if there is a constant $c_1>0$ such that $|s_1-s_2|\geq c_1\min(|s_1-p|,|s_2-p|)$,
then the variation of the measure above is uniformly bounded.
This completely covers the case when $s_1,s_2$ are on different sides of $p$,
but, it  also applies, if   $s_1,s_2$ are on the same sides of $p$, in the appropriate cases.
If  $s_1,s_2$ are on the same sides of $p$ and $\Gamma\subset\Delta$ is a Lipschitz curve
connecting them, then $\dfrac{Z_\Delta(s_1)^{n_1+1}Z_\Delta(s_2)^{n_2+1}}{Z_\Delta(s)^{n_1+n_2+2}}$
is uniformly bounded, as long as
there is a constant $c_2>0$ such that $|s_1-s_2|\leq c_2\min(|s_1-p|,|s_2-p|)$.
Consequently, the variation of
$\dfrac{Z_\Delta(s_1)^{n_1+1}Z_\Delta(s_2)^{n_2+1}}{Z_\Delta(s)^{n_1+n_2+2}}
\Xi_{\Gamma;k}^{(n_1,n_2)}(s)$ is uniformly bounded, either way.

For general Lipschitz curve systems $\Delta$,  due to the properties of $\psi$, as required
in \ref{po:multsys}, the same qualitative picture applies.

(b)
This is just the superposition of the equalities \eqref{colla}, but  we have to use 
the sequential closedness of $A$ in order to interchange integration and operator actions.
\end{proof}
\end{theorem}
It is a rather trivial but important observation that our statements above
can also be applied in the special case when $u=1\in \mathfrak V=\mathbb C$, $A=z\in\mathbb C$;
yielding statements about scalar functions and measures.

\section{Some function theory}
We recall some function theory:
\begin{lemma}\label{lem:boundary}
Suppose that $\gamma:[a,b]\rightarrow\mathbb C$ is a Lipschitz curve, $\Gamma=\gamma([a,b])$.

Suppose that $\mathcal U$ is an open set, $f:\mathcal U\setminus\Gamma\rightarrow\mathbb C$
is holomorphic, and it extends to the analytic boundary $\anb (\mathcal U\setminus\Gamma)|_{\mathcal U}$
continuously.
Then, we claim,
\[-\pi^{-1}\bar\partial f=\gamma_*\nu\]
in distributional sense on $\mathcal U$, with the measure
\[\nu(t)=\frac{f(\gamma(t)^+)-f(\gamma(t)^-)}{2\pi\mathrm i}\,\mathrm d\gamma(t).\]
\begin{proof}
This is a consequence of Green's theorem (cf. \cite{BG},\cite{F}).
\end{proof}
\end{lemma}

\begin{lemma} \label{lem:summary}
Suppose that $\Gamma\subset\mathbb C$ is a compact Lipschitz curve sytem,
$\mu$ is a complex Borel measure of finite variation on $\Gamma$, and
\[f(z)=\int_{s\in\Gamma} \frac{\mu(s)}{s-z}.\]
Then the function $f(z)$ has the following properties:

(a) $f(z)$ is holomorphic on $\mathbb C\setminus\Gamma$, $f(z)$ vanishes at $\infty$.

(b) $|f(z)|\leq\|\mu\|d_{\mathbb C}(z,\Gamma)^{-1}$.
If for some $s\in \Gamma$, $\alpha\in \mathbb R$, $\varepsilon_1,\varepsilon_2>0$
the small sector $\{s+r\mathrm e^{\mathrm i\beta}\,:\, r\in(0,\varepsilon_1), \beta\in(\alpha-\varepsilon_2,
\alpha+\varepsilon_2)\}$ is in $\mathbb C\setminus \Gamma$, then
\[\lim_{t\searrow0} (t\mathrm e^{\mathrm i\alpha}) f(s+t\mathrm e^{\mathrm i\alpha})=\mu(\{s\}).\]

(c) $f(z)$ is a function of class $\mathcal L^p_{\mathrm{loc}}$ for $1\leq p<2$;
$f(z)$ can be understood in distributional sense, and then
\[-\pi^{-1}\bar\partial f=\mu.\]
In particular, if  $f(z)$ is given, then the measure $\mu$ can be recovered.

(d) If $s\in\Gamma$ is a regular point, and in its open neighborhood  $\mathcal U$, the function
$f$ extends to the analytic boundary of $\overline{\mathbb C}\setminus\Gamma$
continuously,  and we orient $\Gamma$ in $\mathcal U$, then on $\mathcal U\cap\Gamma$
\[\mu(s)=\frac{f(s^+)-f(s^-)}{2\pi\mathrm i}\,\mathrm d\Gamma(s).\]
\begin{proof}
(a) is immediate from the definition;
(b) is a simple estimate.
(c): On  $\mathbb C\setminus\Gamma$ the function $f(z)$ is
given as the convolution of the function $z\mapsto-\frac1z$ and the compactly supported
finite measure $\mu$.
The function $z\mapsto-\frac1z$ is of class $\mathcal L^p_{\mathrm{loc}}$ ($1\leq p<2$),
hence we see that  $f(z)$ is a function of class  $\mathcal L^p_{\mathrm{loc}}$ ($1\leq p<2$).
 This convolution can also be understood in  distributional sense.
We know that the function $z\mapsto\frac1{\pi z}$ is a fundamental solution of the $\bar\partial$
operator; from which $-\pi^{-1}\bar\partial f(z)=
\bar \partial*\left(z\mapsto\frac1{\pi z}\right)*\mu=\mu$ follows.
(d) follows from Lemma \ref{lem:boundary} and point (c).
\end{proof}
\end{lemma}

A convenient lemma in the other direction is:
\begin{lemma} \label{lem:reverse}
Suppose that $\gamma:[a,b]\rightarrow\mathbb C$ is a Lipschitz curve, $\Gamma=\gamma([a,b])$.
Assume that $f(z)$ is a function on $\overline{\mathbb C}\setminus\Gamma$ such that

(i) $f$ is holomorphic on $\mathbb C\setminus\Gamma$, and $f$ vanishes at $\infty$;

(ii) $f$ extends to a distribution on $\mathbb C$;

(iii) $f$ extends continuously to the
to the analytic boundaries $\Gamma^+$ and $\Gamma^-$, and  the function
\[m(t)=\frac{f(\gamma(t)^+)-f(\gamma(t)^-)}{2\pi\mathrm i}\]
is in $\mathcal L^1([a,b])$;

(iv) for  s=$\gamma(a)$ and  $s=\gamma(b)$ there exist
$\alpha\in \mathbb R$, $\varepsilon_1,\varepsilon_2>0$ such that
the small sector $\{s+r\mathrm e^{\mathrm i\beta}\,:\, r\in(0,\varepsilon_1), \beta\in(\alpha-\varepsilon_2,
\alpha+\varepsilon_2)\}$ is in $ \mathbb C\setminus \Gamma$, and
\[\lim_{t\searrow0} (t\mathrm e^{\mathrm i\alpha}) f(s+t\mathrm e^{\mathrm i\alpha})=0.\]

Then, we claim,
\begin{equation}f(z)=\int_{t=a}^b m(t)\frac{\mathrm d\gamma(t)}{\gamma(t)-z}.\label{eq:pres}\end{equation}
\begin{proof}
Let, temporarily, $g$ denote the right side of \eqref{eq:pres}.
If $g$ is such a function, then, due to Lemma \ref{lem:summary}.d,   the distribution $\bar\partial(g-f)$ is a
distribution supported at the endpoints of $\gamma$.
Then consider $r=\left(z\mapsto\frac1{\pi z}\right)*\bar\partial(g-f)$.
We see that $r$ is a rational function with singular support at the endpoints of $\gamma$ at most,
and one which vanishes at  $\infty$.
Also, $\bar\partial$ vanishes on $(g-f)-r$.
Then, from Liouville's theorem, we see that $g-f=r$.
Due to Lemma \ref{lem:summary}.b and condition (d) here, we see that $g-f$ has actually no singular support,
consequently $g-f=0$.
\end{proof}
\end{lemma}

\section{The holomorphic functional calculus}

\begin{point}[The idea of the functional calculus]
Suppose that $\bo\mu$ is a resolvent measure on a Lipschitz curve sytem $\Gamma$.
We can substitute any $z\in\overline{\mathbb C}\setminus\Gamma$ to the place of $X$,
and we get a well-defined function
\[f(z)=\int_{s\in\Gamma} \bo\mu(s,z).\]
We call such functions slowly growing holomorphic functions on $\overline{\mathbb C}\setminus\Gamma$
of order $n$, and we  denote their set  by $\mathcal H^n_\Gamma$.
We must notice that the same function might belong to  very different resolvent measures.

More generally, if $u\in\mathrm D_{\Gamma}^n(A)$, then we can take
\[f(A)u\stackrel{?}{=}\int_{s\in\Gamma} \bo\mu(s,A)u.\]
Indeed, it is natural for this integral to be considered as $f(A)u$, except we have to show that
this expression depends only on $f$ but not on the resolvent measure used.
\end{point}
\begin{point}[Global conformal change] \label{po:confcha}
If $h(z)=\frac{az+b}{cz+d}$, $ad-bc\neq 0$ is a global conformal map, then
one can define a resolvent measure $h_*\bo\mu$ on $h(\Gamma)$
by $h_*\bo\mu(t,Y)=\mu(h^{-1}(t),h^{-1}(X))$.  Then
\[\int_{t\in h(\Gamma)} (h_*\bo\mu)(t,h_{\mathrm f}(A))u=\int_{s\in\Gamma} \bo\mu(s,A)u.\]
If $\overline\Gamma, h(\overline\Gamma)\subset\mathbb C$, then this yields
\[h_*\bo\mu(t,Y)=\sum_{k=0}^n\left(\frac{c}{c\,h^{-1}(t)+d}+\frac{ad-bc}{(c\,h^{-1}(t)+d)^2}
\frac{1}{t-Y}\right)^k\mu_k(h^{-1}(t)),\]
cf. Lemma \ref{qlem:confchange}; but it can also be done in general.

This implies, however, that, practically, we can always assume that $\overline\Gamma\subset\mathbb C$,
because we can always transform $\overline\Gamma$ into such one by a  conformal transformation $h$.
Then we must pass from $A$ to $h_{\mathrm f}(A)$, and instead of $f(z)$, we must deal with
\[f\circ h^{-1}(z)= \int_{t\in h(\Gamma)} h_*\bo\mu(t,z).\]
(And this  yields $f\circ h^{-1}(h_{\mathrm f}(A))u=f(A)u$ once well-definedness is established.)
\end{point}
So, in what follows, it is sufficient to deal with the case $\overline{\Gamma}\subset\mathbb C$.
\begin{point}
Suppose that $\Gamma$ is a Lipschitz curve system,
$\overline\Gamma\subset\mathbb C$,
with connected components $\Gamma_i$, $i=1,\ldots,b_0(\Gamma)$, $w_i\in\Gamma_i$, $n\geq 1$.
We say that $\bo\mu$ is a resolvent measure of normal form of order $n$ with respect to the $w_i$, if it is of shape
\begin{equation}
c_0+\sum_{i=1}^{b_0(\Gamma)}\sum_{j=1}^{n-1}c_{w_i,j}\frac{\bo\delta_{w_i}(s)}{(s-X)^j}
+\frac{\mu_n(s)}{(s-X)^n}.\label{eq:normalform}
\end{equation}
Strictly speaking, this is not a resolvent measure, as we have not specified
where the scalar measure is supported, but it causes no problems when we integrate the resolvent measure.
Then its complex evaluation
on $\overline{\mathbb C}\setminus\Gamma $ is the function
\[f(z)=c_0+\sum_{i=1}^{b_0(\Gamma)}\sum_{j=1}^{n-1}\frac{c_{w_i,j}}{(w_{i}-z)^j}
+\int_{s\in\Gamma}\frac{\mu_n(s)}{(s-z)^n}.\]

Suppose $\mathrm d\Delta_k$, $k=1,\ldots,b_1(\Gamma)$ are arc variation measures of
a representative set of independent (oriented) loops.
We say that the resolvent measure $\bo\mu$ is a circular measure of order $n$, if it is
of shape
\[\sum_{k=1}^{b_1(\Gamma)}p_i(s)\frac{\mathrm d\Delta_k(s)}{(s-X)^n},\]
where $p_i$ are complex polynomials of order less than $n-1$.
The measures  $\mathrm d\Delta_k$ can, of course, be replaced by the arc variation measures
 $\mathrm d\anb \mathcal U_k$ of the bounded connected components of $\overline{\mathbb C}\setminus\Gamma$.
\end{point}
\begin{lemma}\label{lem:normalf}
If $\bo\mu$ is a resolvent measure on $\Gamma$ of normal form of degree $n$, and
$f(z)=\int_{s\in\Gamma} \bo\mu(s,z)$ is identically zero on $\overline{\mathbb C}\setminus\Gamma$ , then
$\bo\mu$ is a circular measure of order $n$.
\begin{proof} First we prove the case when $\Gamma$ is compact.
Assume that $\bo\mu$ is of shape \eqref{eq:normalform}.
As $f$ vanishes in $\infty$, we know that $c=0$.
Taking various line integrals around the connected components of $\Gamma$, we can deduce
that the residues vanish, so $c_{w_i,1}=0$.
This implies that $f(z)$ allows the  primitive function
\[ f^{[[1]]}(z)=\sum_{i=1}^m\sum_{j=2}^{n-1}(j-1)\frac{c_{w_i,j}}{(w_{i}-z)^{j-1}}
+\int_{s\in\Gamma}(n-1)\frac{\mu_n(s)}{(s-z)^{n-1}}. \]
This gives  constant functions on the connected components of  $\overline{\mathbb C}\setminus\Gamma$
but still vanishes on the component containing $\infty$.
Iterating this argument, we find that all $c_{w_i,j}$=0, and $f(z)$ allows the $(n-1)$th primitive function
\[f^{[[n-1]]}(z)=\int_{s\in\Gamma}(n-1)!\frac{\mu_n(s)}{s-z}.\]
This function vanishes on the unbounded component of $\overline{\mathbb C}\setminus\Gamma$,
and it gives polynomials $q_k$ of degree less than $n-1$ on the bounded components $\mathcal U_k$, $k=1,\ldots,b_1(\Gamma)$. Then, by Lemma \ref{lem:summary}.d, we know its behaviour except at the nodes:
\[(n-1)!\mu_n(s)= \sum_{k=1}^{b_1(\Gamma)}\frac{q_k(s)}{2\pi\mathrm i}\mathrm d\anb \mathcal U_k
+\sum_{p\in N(\Gamma)}\tilde c_p\bo\delta_p(s).\]
But the the polynomial behaviour of $f$ forces the $\tilde c_p$ to vanish, so
\[\bo\mu(s,X)=\frac{\mu_n(s)}{(s-X)^n}= \sum_{k=1}^{b_1(\Gamma)}\frac1{(n-1)!}
\frac{q_k(s)}{2\pi\mathrm i}\frac{\mathrm d\anb \mathcal U_k(s)}{(s-X)^n}.\]

If $\Gamma$ is non-compact, then, applying Lemma \ref{lem:qRIF}.b, we can prepare a normal form $\bo\mu'$
on $\overline\Gamma$ with corrections which are finitely many polynomial times arc variation measures, apart from
finitely many singular points from $\Gamma$. This normal form $\bo\mu'$ on $\overline\Gamma$
must be circular; in particular $\bo\mu$ must be polynomial times arc variation measure.
But this polynomiality implies that $\mu$ must vanish around the points of $\overline\Gamma\setminus\Gamma$.
In particular, $\mu$ descends to $\Gamma_{\mathrm{amp}}$, where some short pieces of the legs
joining to points of   $\overline\Gamma\setminus\Gamma$ are amputated.
$\Gamma_{\mathrm{amp}}$ is compact, so $\bo\mu$ must be circular.
\end{proof}
\end{lemma}
Now, we are in position to prove our main theorem:

\begin{theorem}
Suppose that $\mathfrak V$ is a sequentially complete locally convex vector space, and $A$
is a multivalued linear operator on $\mathfrak V$.

(a)
If $f\in\mathcal H_\Gamma^n$, $u\in\mathrm D^n_\Gamma(A)$, then
$f(A)u$ is well-defined, i.~e.~any way we choose a resolvent measure $\bo\mu$ of order
$n$ on $\Gamma$ such that
\[f(z)=\int_{s\in\Gamma}\bo\mu(s,z),\]
the computed value
\[f(A)u=\int_{s\in\Gamma}\bo\mu(s,A)u\]
does not depend on $\bo\mu$.
This construction is linear in $f$ and $u$.

(b) If $f\in\mathcal H_\Gamma^n$, $u\in\mathrm D^n_\Gamma(A)$,
and $h:\overline{\mathbb C}\rightarrow\overline{\mathbb C}$ is a global conformal map, then
$f\circ h^{-1}\in\mathcal H_{h(\Gamma)}^n$, $u\in\mathrm D^n_{h(\Gamma)}(h_{\mathrm f}(A))$, and
\[f(A)u=(f\circ h^{-1})(h_{\mathrm f}(A))u.\]

(c) If $f\in\mathcal H_\Gamma^n$, and $A$ is sequentially closed, then $f(A)$ defines continuous linear maps
\[f(A): \mathrm D_{\Gamma}^{n+n_1}(A)\rightarrow \mathrm D_{\Gamma}^{n_1}(A).\]

(d) If $A$ is a sequentially closed,
and $f_1\in\mathcal H_\Gamma^{n_1}$, $f_2\in\mathcal H_\Gamma^{n_2}$,
then $f_1f_2\in\mathcal H_\Gamma^{n_1+n_2}$, and for any $u\in\mathrm D_{\Gamma}^{n_1+n_2}(A)$,
\[(f_1f_2)(A)u=f_1(A)f_2(A)u.\]

\begin{proof}
(a) It is sufficient to prove that if $\bo\mu$ gives $f(z)\equiv0$, then it also
yields $f(A)u=0$. By a repeated use of Theorem  \ref{thm:addred},
we can assume that  $\bo\mu$ is of  normal form.
But then, by Lemma \ref{lem:normalf}, $\bo\mu$ must be circular.
Then, however, by Lemma \ref{lem:qRIF}.c,  $f(A)u=0$.

(b) This follows from the discussion in \ref{po:confcha}.

(c) This follows from the discussion in \ref{po:measint}.

(d) This follows from Theorem \ref{thm:multred}.
\end{proof}
\end{theorem}
\begin{point} We see that $\mathcal H_\Gamma=\bigcup_{n\in\mathbb N}\mathcal H_\Gamma^n$ forms
a filtered algebra. The elements of $\mathcal H_\Gamma$ can be extended from $\overline{\mathbb C}\setminus\Gamma$
to $\overline{\mathbb C}$ as distributions. Indeed, if $\overline\Gamma\subset\mathbb C$, what we can practically
assume, then such an extension
is given by
$\left(\int_{s\in\Gamma}\mu_0(s)\right)+\left(z\rightarrow -\frac1z\right)*\sum_{k=1}^n\frac{\partial^{k-1}}{(k-1)!}\mu_k$
with respect to \eqref{eq:resmea}.
\end{point}

\section{Example: Logarithms and fractional powers along Lipschitz curves}
\begin{point}
Let $\gamma:[a,b]\rightarrow\mathbb C$ be a Lipschitz   curve, $\Gamma=\gamma([a,b])$.
Then, for $z\in\overline{\mathbb C}\setminus\Gamma$, we define the associated logarithm function as
\[\log(\gamma,z)=\int_{t=a}^b \frac{\mathrm d\gamma(t)}{\gamma(t)-z}
=\log|\gamma(b)-z|-\log|\gamma(a)-z|+\mathrm i\int_{t=a}^b \mathrm d\Arg(\gamma(t)-z).\]
In fact, $\log(\gamma,z)=\widetilde\log\dfrac{\gamma(b)-z}{\gamma(a)-z}$,
where $\widetilde\log$ is the unique branch of logarithm such that the function vanishes at $\infty$.
Regarding local behaviour, for $z\sim\gamma(s)$, $s\in(a,b)$ we see that
\begin{equation}
\log(\gamma,z)=\widehat\log\dfrac{\gamma(b)-z}{z-\gamma(a)}+\sgn(\gamma,z)\pi\mathrm i,
\label{eq:log}
\end{equation}
where $\widehat\log$ is an appropriate branch of $\log$, and $\sgn(\gamma,z)$ is $+1$ or $-1$,
if $z$ is on the  left or right side of $\gamma$, respectively.
~\\

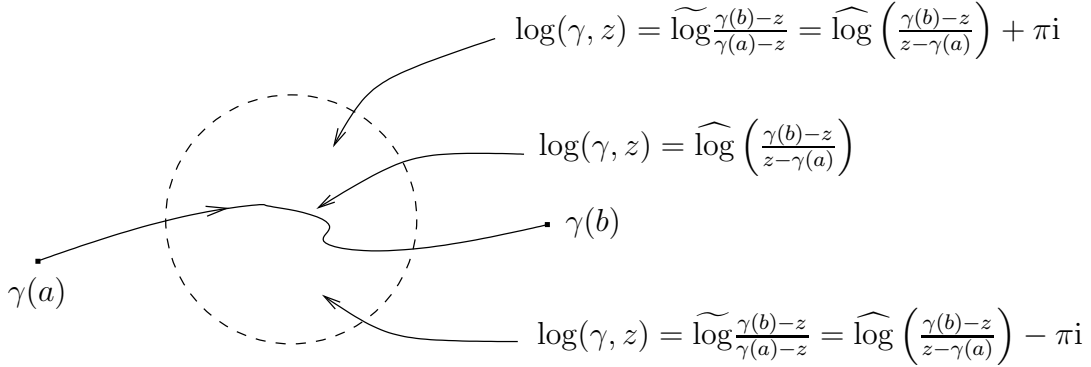
\input{el4.pstex_t}
~\\

Considering \eqref{eq:log}, it seems reasonable to extend the logarithm function to
$z=\gamma(s)\in\gamma((a,b))$ by taking $\sgn(\gamma,z)=0$.

In order to avoid technical difficulties, we simply assume that $\gamma$ is differentiable at
its endpoints $\gamma(a)$, $\gamma(b)$. This makes
$\log(\gamma,z)=\log^+\left(\frac1{|\gamma(a)-z|}\right)-\log^+\left(\frac1{|\gamma(b)-z|}\right)+
\textrm{a bounded term}$.
\end{point}

\begin{point} \label{po:logref1}
(a) From Lemma \ref{lem:reverse}, we can deduce that for $n\in\mathbb N$,
\[\log(\gamma,z)^{n+1}=\int_{t=a}^b \frac1{2\pi\mathrm i}
\left(  (\log(\gamma,\gamma(t)) +\pi\mathrm i)^{n+1}-
(\log(\gamma,\gamma(t)) -\pi\mathrm i)^{n+1} \right) \frac{\mathrm d\gamma(t)}{\gamma(t)-z}.\]
In particular, $\log(\gamma,z)^{n+1}\in\mathcal H_{\Gamma}^1$.

(b) Then, by simple summation, we find that for $|\alpha|<1$,
\[z^\alpha_\gamma\equiv\exp \left(\alpha\log(\gamma,z)\right)=1+\int_{t=a}^b\frac{\sin \alpha\pi }{\pi}
\exp \left(\alpha\log(\gamma,\gamma(t))\right)\frac{\mathrm d\gamma(t)}{\gamma(t)-z}.\]
In particular, $\exp \left(\alpha\log(\gamma,z)\right)\in\mathcal H_{\Gamma}^1$.

In fact, it is transparent that this latter formula extends to $-1<\Rea \alpha<1$ by analytic continuation, and
$z^\alpha_\gamma\in\mathcal H_{\gamma([a,b))}$ for $0<\Rea \alpha<1$, and
$z^\alpha_\gamma\in\mathcal H_{\gamma((a,b])}$ for $-1<\Rea \alpha<0$.
\end{point}

\begin{point}\label{po:logref2}
In the special case when $[a,b]=[-1,1]$ and $\gamma(t)=t$, it yields
\begin{multline}
\left(\log\left(\frac{z-1}{z+1}\right)\right)^{n+1}=\\=\int_{t=-1}^1 \frac1{2\pi\mathrm i}
\left(\left(\log\left(\frac{1-t}{1+t}\right) +\pi\mathrm i\right)^{n+1}-
\left(\log\left(\frac{1-t}{1+t}\right) -\pi\mathrm i\right)^{n+1}\right) \frac{\mathrm dt}{t-z},
\notag\end{multline}
and
\[\left(\frac{z-1}{z+1}\right)^{\alpha}
=1+\int_{t=-1}^1\frac{\sin \alpha\pi }{\pi}\left(\frac{1-t}{1+t}\right)^{\alpha}
\frac{\mathrm dt}{t-z}.\]

We can apply a global conformal transformation, $w=\frac{z-1}{z+1}$, $r=\frac{t-1}{t+1}$;
which yields
\[(\log w)^{n+1}=\int^{0}_{r=-\infty}\frac{(\log( -r)+\pi\mathrm i)^{n+1}-(\log( -r)-\pi\mathrm i)^{n+1}}{2\pi\mathrm i}\frac{w-1}{(1-r)(w-r)}\mathrm dr,\]
and
\[w^\alpha=1+\int^{0}_{r=-\infty}\frac{\sin\alpha\pi}\pi(-r)^\alpha\frac{w-1}{(1-r)(w-r)}\mathrm dr.\]
(Here, from the viewpoint of the functions,  we cut the complex plane along the negative real axis.)
The same qualitative behaviour remains:  $z^\alpha\in\mathcal H_{[-\infty,0)}^1$ for
$0<\Rea \alpha<1$, etc.

In a sense, the ordinary $\log$ function is the special case  when $\gamma:[a,b]\rightarrow\overline{\mathbb C}$
is a Lipschitz curve from $-\infty$ to $0$ along the negative real axis.
Strictly,  this is true only up to an additive constant.
$\log(\gamma,z)$ is characterized by the orientation of $\Gamma$, the fact that the difference
between the left and right boundary values of $\Gamma$ is exactly $2\pi\mathrm i$, and that
$\log(\gamma,\infty)=0$. This latter condition is not invariant under global conformal changes,
but it is only a matter of an additive constant. Anyway, it is safe to associate
this case to the straight Lipschitz curve from $-\infty$ to $0$. This illustrates the general case:
\end{point}
\begin{point}
If $\gamma:[a,b]\rightarrow\overline{\mathbb C}$ is a Lipschitz curve, then we can define
$\log(\gamma,z)$ on $ \overline{\mathbb C}\setminus\gamma([a,b])$, although the
definition is not entirely canonical, but up to an additive constant. (So, the logarithm
is not necessarily ``natural''.) Then $z^\alpha_\gamma\equiv\exp \left(\alpha\log(\gamma,z)\right)$
can also be considered.
\end{point}

\begin{theorem} Suppose that $A$ is a sequentially closed multivalued operator on the
sequentially complete locally convex vector space $\mathfrak V$, and
 $\gamma:[a,b]\rightarrow\overline{\mathbb C}$ is a Lipschitz curve with differentiable endpoints.

(a) If $u\in\mathrm D^2_{\gamma([a,b])}(A)$, then $\log(\gamma,A)$ can be iterated on $u$ arbitrarily many times.
Formal analogues like
\[(\log A)^{n+1}=\int^{0}_{r=-\infty}\frac{(\log( -r)+\pi\mathrm i)^{n+1}-(\log( -r)-\pi\mathrm i)^{n+1}}{2\pi\mathrm i}\frac{A-1}{(1-r)(A-r)}u\, \mathrm dr,\]
($n\in\mathbb N$) hold.

(b)
Suppose that  $n\in\mathbb N$, and $S_n=[-n,n]$ or $S_n=\{\alpha\,:\,|\Rea \alpha|<n\}$.
Then we claim: If  $u\in\mathrm D^n_{\gamma([a,b])}(A)$, then the function
\[\alpha\in S\mapsto A^\alpha_\gamma u\in\mathfrak V\]
is well-defined and continuous.
$A_\gamma^{\alpha_1}A_\gamma^{\alpha_2} u=A_\gamma^{\alpha_1+\alpha_2}u$,
as long as $\alpha_1\in S_{n_1}$, $\alpha_2\in S_{n_2}$, $n_1+n_2\leq n$.
If $\Rea\alpha,\Rea \alpha_i> 0$, then  it is sufficient to assume that
$u\in\mathrm D^n_{\gamma([a,b))}(A)$.

Formal analogues like
\[A^\alpha u=u+\int^{0}_{r=-\infty}\frac{\sin\alpha\pi}\pi(-r)^\alpha\frac{A-1}{(1-r)(A-r)}u\,\mathrm dr,\]
($-1<\Rea \alpha<1$) hold.

(c) Suppose that $\tilde S_n=\{\alpha\,:\,|\alpha|<n\}$, $n\geq 2$, $u\in\mathrm D^n_{\gamma([a,b])}(A)$.
Then for $\alpha\in\tilde S_n$,
\[A^\alpha_\gamma u=\sum_{n=0}^\infty\frac1{n!}\log(\gamma,A)^nu\]
locally uniformly.

(Remark:
$u\in\mathrm D^1_{[-\infty,0]}(A)$ means that  $\frac{(1-r)(A-1)}{A-r}u$ is bounded and continuous,
$r\in[-\infty,0]$; and
$u\in\mathrm D^1_{[-\infty,0)}(A)$ means that  $\frac{(-r)(A-1)}{A-r}u$ is bounded and continuous,
$r\in[-\infty,0)$.)

\begin{proof}
By a global conformal change, we can assume that $\gamma([a,b])\subset\mathbb C$,
where the statement is transparent according the previous discussion and the multiplicative properties
of the functional calculus, except (c).

(c):
Let $h_m(t)=\frac1{2\pi\mathrm i}\left(
(\log(\gamma,\gamma(t))+\pi\mathrm i)^{m+1}-(\log(\gamma,\gamma(t))-\pi\mathrm i)^{m+1}\right)$.
Then
\[A^\alpha_\gamma u=A^{\alpha/n}_\gamma\ldots A^{\alpha/n}_\gamma u=\]
\[=\int_{t_1=a}^{b}\ldots\int_{t_n=a}^{b}
\frac{(\alpha/n)^{m_1}}{m_1!}\ldots\frac{(\alpha/n)^{m_n}}{m_n!}h_{m_1}(t)\ldots h_{m_n}(t)
\frac{\mathrm d\gamma(t_n)}{\gamma(t_n)-A}\ldots\frac{\mathrm d\gamma(t_1)}{\gamma(t_1)-A}u=\]
\[=\sum_{m_1,\ldots,m_n=0}^\infty \frac{(\alpha/n){m_1}}{m_1!}\ldots\frac{(\alpha/n)^{m_n}}{m_n!}
(\log(\gamma,A))^{m_1}\ldots(\log(\gamma,A))^{m_n}u=\]
\[=\sum_{m=0}^\infty \frac{\alpha^m}{m!}(\log(\gamma,A))^mu.\]
\end{proof}
\end{theorem}

\end{document}

%% file: el4.pstex_t
\begin{picture}(0,0)%
\includegraphics{el4.pstex}%
\end{picture}%
\setlength{\unitlength}{4144sp}%
\begingroup\makeatletter\ifx\SetFigFont\undefined%
\gdef\SetFigFont#1#2#3#4#5{%
  \reset@font\fontsize{#1}{#2pt}%
  \fontfamily{#3}\fontseries{#4}\fontshape{#5}%
  \selectfont}%
\fi\endgroup%
\begin{picture}(3328,2050)(221,-1888)
\put(236,-1559){\makebox(0,0)[lb]{\smash{{\SetFigFont{12}{14.4}{\rmdefault}{\mddefault}{\updefault}{\color[rgb]{0,0,0}$\gamma(a)$}%
}}}}
\put(3534,-1144){\makebox(0,0)[lb]{\smash{{\SetFigFont{12}{14.4}{\rmdefault}{\mddefault}{\updefault}{\color[rgb]{0,0,0}$\gamma(b)$}%
}}}}
\put(3379,-689){\makebox(0,0)[lb]{\smash{{\SetFigFont{12}{14.4}{\rmdefault}{\mddefault}{\updefault}{\color[rgb]{0,0,0}$\log(\gamma,z)=\widehat\log\left(\frac{\gamma(b)-z}{z-\gamma(a)}\right)$}%
}}}}
\put(3368,-1819){\makebox(0,0)[lb]{\smash{{\SetFigFont{12}{14.4}{\rmdefault}{\mddefault}{\updefault}{\color[rgb]{0,0,0}$\log(\gamma,z)=\widetilde\log\frac{\gamma(b)-z}{\gamma(a)-z}=\widehat\log\left(\frac{\gamma(b)-z}{z-\gamma(a)}\right)-\pi\mathrm i$}%
}}}}
\put(3238,  3){\makebox(0,0)[lb]{\smash{{\SetFigFont{12}{14.4}{\rmdefault}{\mddefault}{\updefault}{\color[rgb]{0,0,0}$\log(\gamma,z)=\widetilde\log\frac{\gamma(b)-z}{\gamma(a)-z}=\widehat\log\left(\frac{\gamma(b)-z}{z-\gamma(a)}\right)+\pi\mathrm i$}%
}}}}
\end{picture}%